%% file: Hamkins_main.tex
\documentclass[12pt,runningheads]{fotfs}

\usepackage{amsmath, amssymb, graphics, latexsym, makeidx}
\usepackage{times}

\makeindex
\begin{document}
\frontmatter
\mainmatter
\include{Hamkins_paper}
\backmatter
\newpage
\pagestyle{empty}
\printindex
\newpage
\pagestyle{empty}
\include{Hamkins_bib}

\end{document}

%% file: Hamkins_paper.tex

\def\col#1#2#3{\hbox{\vbox{\baselineskip=0pt\parskip=0pt\cell#1\cell#2\cell#3}}}
\newcommand{\cell}[1]{\boxit{\hbox to 17pt{\strut\hfil$#1$\hfil}}}
\newcommand{\head}[2]{\lower2pt\vbox{\hbox{\strut\footnotesize\it\hskip3pt#2}\boxit{\cell#1}}}
\newcommand{\boxit}[1]{\setbox4=\hbox{\kern2pt#1\kern2pt}\hbox{\vrule\vbox{\hrule\kern2pt\box4\kern2pt\hrule}\vrule}}
\newcommand{\Col}[3]{\hbox{\vbox{\baselineskip=0pt\parskip=0pt\cell#1\cell#2\cell#3}}}
\newcommand{\tapenames}{\raise 5pt\vbox to .6in{\hbox to .8in{\it\hfill input: \strut}\vfill\hbox to
.8in{\it\hfill scratch: \strut}\vfill\hbox to .8in{\it\hfill output: \strut}}}
\newcommand{\Head}[4]{\lower2pt\vbox{\hbox to25pt{\strut\footnotesize\it\hfill#4\hfill}\boxit{\Col#1#2#3}}}
\newcommand{\Dots}{\raise 5pt\vbox to .6in{\hbox{\ $\cdots$\strut}\vfill\hbox{\ $\cdots$\strut}\vfill\hbox{\
$\cdots$\strut}}}
\renewcommand{\dots}{\raise5pt\hbox{\ $\cdots$}}
\newcommand{\df}{\it}
\newcommand\R{{\mathbb R}}
\newcommand{\N}{{\mathbb N}}
\newcommand{\from}{\mathbin{\vbox{\baselineskip=3pt\lineskiplimit=0pt
                         \hbox{.}\hbox{.}\hbox{.}}}}
\newcommand{\of}{\subseteq}
\newcommand{\restrict}{\upharpoonright}
\def\Theorem #1.#2 #3\par{\setbox1=\hbox{#1}\ifdim\wd1=0pt
   \begin{theorem}{\rm #2} #3\end{theorem}\else
   \newtheorem{#1}[theorem]{#1}\begin{#1}\label{#1}{\rm #2} #3\end{#1}\fi}
\def\Corollary #1.#2 #3\par{\setbox1=\hbox{#1}\ifdim\wd1=0pt
   \begin{corollary}{\rm #2} #3\end{corollary}\else
   \newtheorem{#1}[theorem]{#1}\begin{#1}\label{#1}{\rm #2} #3\end{#1}\fi}
\def\Lemma #1.#2 #3\par{\setbox1=\hbox{#1}\ifdim\wd1=0pt
   \begin{lemma}{\rm #2} #3\end{lemma}\else
   \newtheorem{#1}[theorem]{#1}\begin{#1}\label{#1}{\rm #2} #3\end{#1}\fi}
\def\Question #1.#2 #3\par{\setbox1=\hbox{#1}\ifdim\wd1=0pt
   \begin{question}{\rm #2} #3\end{question}\else
   \newtheorem{#1}[theorem]{#1}\begin{#1}\label{#1}{\rm #2} #3\end{#1}\fi}
\def\Proof: {\begin{proof}}
\newcommand{\QED}{\end{proof}}
\newcommand{\WO}{\mathop{\hbox{\sc wo}}}
\newcommand{\jump}{{\!\triangledown}}
\newcommand{\Jump}{{\!\blacktriangledown}}
\def\ilt{<_{\infty}}
\def\ileq{\leq_{\infty}}
\def\iequiv{\equiv_{\infty}}
\newcommand{\dom}{\mathop{\rm dom}}
\def\<#1>{\langle\,#1\,\rangle}
\newcommand{\set}[1]{\{\,{#1}\,\}}
\newcommand{\st}{\mid}

\title{Supertask Computation}

\author{Joel David Hamkins}

\institute{The College of Staten Island of The City University of New York\\
           Mathematics Department 1S-215, 2800 Victory Boulevard, Staten Island, NY 10314 USA\\
           and\\
           The Graduate Center of The City University of New York\\
           Mathematics Program, 365 Fifth Avenue, New York, NY 10016 USA\\
           http://jdh.hamkins.org}

\email{jdh@hamkins.org}

\grant{My research has been supported in part by grants from the PSC-CUNY Research Foundation and from the
NSF.}

\received{...}
\revised{...}
\accepted{...}

\subjclass{{\bf 03D10} 68Q05 03D25}

\maketitle
\begin{abstract}
Infinite time Turing machines extend the classical Turing machine
concept to transfinite ordinal time, thereby providing a natural model
of infinitary computability that sheds light on the power and
limitations of supertask algorithms.
\end{abstract}

\section{Supertasks}

What would you compute with an infinitely fast computer? What {\it could} you compute? To make sense of
these questions, one would want to understand the algorithms that the machines would carry out,
computational tasks involving infinitely many steps of computation. Such tasks, known as {\it supertasks},
have been studied since antiquity from a variety of viewpoints.

Zeno of Elea (ca. 450 B.C.) was perhaps the first to grapple with supertasks, in his famous paradox that it
is impossible to go from here to there, because before doing so one must first get halfway there, and
before that halfway to the halfway point, and so on, {\it ad infinitum}. Zeno takes the impossibility of
completing a supertask as the foundation of his reductio. More recently, twentieth century philosophers
(see \cite{Hamkins_Thomson}) have introduced Thomson's lamp, which is on for 1/2 minute, off for 1/4
minute, on for 1/8 minute, and so on. After one minute, is it on or off?

In a more intriguing example, let's suppose that you have infinitely many one dollar bills (numbered $1$,
$3$, $5$, $\cdots$) and in some nefarious underground bar, the Devil explains to you that he has an
attachment to your particular bills, and is willing to pay you two dollars for each of your one dollar
bills. To carry out the exchange, he proposes an infinite series of transactions, in each of which he will
hand over to you two dollars and take from you one dollar. The first transaction will take 1/2 hour, the
second 1/4 hour, the third 1/8 hour, and so on, so that after one hour the entire exchange is complete.
Should you accept his proposal? Perhaps you will become richer? At the very least, you think, it will do no
harm, and so the contract is signed and the procedure begins.

It appears initially that you have made a good bargain, because at
every step of the transaction, you receive two dollars but give up only
one. The Devil is particular, however, about the order in which the
bills are exchanged: he always buys from you your lowest-numbered bill,
paying you with higher-numbered bills. (So on the first transaction he
accepts from you bill number 1, and pays you with bills numbered 2 and
4, and on the second transaction he buys from you bill number 2, which
he had just paid you, and pays you bills numbered 6 and 8, and so on.)
When the transaction is complete, you discover that you have no money
left at all! The reason is that at the $n^{th}$ exchange, the Devil
took from you bill number $n$, and never subsequently returned it to
you. Thus, the final destination of every individual bill is under the
ownership of that shrewd banker, the Devil.

The point is that you should have paid more attention to the details of the supertask transaction that you
had agreed to undertake. And similarly, when we design supertask algorithms to solve mathematical
questions, we must take care not to make inadvertent assumptions about what may be true only for finite
algorithms.

Supertasks have also been studied by the physicists (see \cite{Hamkins_Earman95}). Using only the Newtonian
gravity law (and neglecting relativity), it is possible to arrange finitely many stars in orbiting pairs,
each pair orbiting the common center of mass of all the pairs, and a single tiny moon racing faster around,
squeezing just so between the dual stars so as to pick up speed with every such transaction. Assuming point
masses (or collapsing stars to avoid collision), the arrangement leads by Newton's law of gravitation to
infinite acceleration in finite time. Other supertasks reveal apparent violations of the conservation of
energy in Newtonian physics: infinitely many billiard balls, of successively diminishing size converging to
a point, are initially at rest, but then the first is set rolling, and each ball transfers in turn all the
energy to the next; after a finite amount of time, all motion has ceased, though every interaction is
energy conserving. Still other arrangements have the balls spaced out further and further out to infinity,
and the interesting thing about both of these examples is that time-symmetry allows them to run in reverse,
with static configurations of balls suddenly coming into motion without violating conservation of energy in
any interaction.

More computationally significant supertasks have been proposed by
physicists in the context of relativity theory
(\cite{Hamkins_EarmanNorton93}, \cite{Hamkins_Hogarth92},
\cite{Hamkins_Hogarth94}). Suppose that you want to know the answer to
some number theoretic conjecture, such as whether there are additional
Fermat primes (primes of the form $2^{2^n}+1$), a conjecture that can
be confirmed with a single numerical example. The way to solve the
problem is to board a rocket, while setting your graduate students to
work on earth looking for an example. While you fly faster and faster
around the earth, your graduate students, and their graduate students
and so on, continue the exhaustive search, with the agreement that if
they ever find an example, they will send a radio signal up to the
rocket. The point is that meanwhile, by accelerating sufficiently fast
towards the speed of light, it is possible to arrange that because of
relativistic time contraction, what is a finite amount of time on the
rocket corresponds to an infinite amount of time on the earth. The
general observation is that by means of such communication between two
reference frames, what corresponds to an infinite search can be
completed in a finite amount of time.

Even more complicated arrangements, with rockets flying around rockets,
can be arranged to solve more complicated number theoretic questions.
And more complicated relativistic spacetimes can be (mathematically)
constructed to avoid the unpleasantness of infinite acceleration
required in the rocket examples above (see \cite{Hamkins_Pitowsky}).

These computational examples speak to Church's thesis, the widely accepted philosophical principle that the
classical theory of computability has correctly captured the notion of what it means to be computable.
Because the relativistic rocket examples provide algorithms for computing functions, such as the halting
problem, that are not computable by Turing machines, one can view them as refuting Church's thesis.
Supporters of this view emphasize that when thinking about what is in principle computable, we must attend
to the computational power available to us as a consequence of the fact that we live in a relativistic or
quantum-mechanical universe. To ignore this power is to pretend that we live in a Newtonian world. Another
simpler argument against Church's thesis consists of the observation that a particle undergoing Brownian
motion can be used to generate a random bit stream that we have no reason to think is recursive. Therefore,
proponents argue, we have no reason to believe Church's thesis.

Apart from the question of what one can actually compute in this world,
whether Newtonian or relativistic or quantum-mechanical, mathematicians
are interested in what {\it in principle} a supertask can accomplish.
Buchi \cite{Hamkins_Buchi62} and others initiated the study of
$\omega$-automata and Buchi machines, involving automata and Turing
machine computations of length $\omega$ which accept or reject infinite
input. Moving to a higher level in the hierarchy, Gerald Sacks and many
others (see \cite{Hamkins_Sacks}) founded the field of higher recursion
theory, including $\alpha$-recursion and $E$-recursion, a huge body of
work analyzing computation on infinite objects. Blum, Shub and Smale
\cite{Hamkins_BlumShubSmale} have presented a model of computation on
the real numbers, a kind of flowchart machine where the basic units of
computation consist of real numbers, in full glorious precision. Apart
from this previous mathematical work, I would like to propose here a
new model of infinitary computability: infinite time Turing machines.
This model offers the strong computational power of higher recursion
theory while remaining very close in spirit to the computability
concept of ordinary Turing machines.

\section{Infinite time Turing machines}

I propose to extend the Turing machine concept to transfinite ordinal time, thereby providing a natural
model for infinitary computability.\footnote{Infinite time Turing machines were originally defined by Jeff
Kidder in 1990, and he and I worked out the early theory together while we were graduate students at UC
Berkeley. Later, Andy Lewis and I solved some of the early questions, and presented a complete introduction
in \cite{Hamkins_InfiniteTimeTM}, later solving Post's problem for supertasks in
\cite{Hamkins_PostProblem}. Benedikt Loewe \cite{Hamkins_Loewe2001}, Dan Seabold \cite{Hamkins_OneTape} and
especially Philip Welch \cite{Hamkins_Welch99}, \cite{Hamkins_Welch2000}, \cite{Hamkins_WelchEventual} have
also made important contributions.} The idea is to allow somehow a Turing machine to compute for infinitely
many steps, while preserving the information produced up to that point.

So let me explain specifically how the machines work. The machine hardware is identical to a classical
Turing machine, with a head moving back and forth reading and writing zeros and ones on a tape according to
the rigid instructions of a finite program, with finitely many states. What is new is the transfinite
behavior of the machine, behavior providing a natural theory of computation on the reals that directly
generalizes the classical finite theory to the transfinite.
\begin{figure}[h]
$$\tapenames\Head100{start}\Col100\Col000\Col100\Col100\Col000\Dots$$
\caption{An infinite time Turing machine: the computation begins}
\end{figure}
For convenience, the machines have three tapes---one for the input, one for scratch work and one for the
output---and the computation begins with the input written out on the input tape, with the head on the
left-most cell in the {\it start} state. The successor steps of computation proceed in exactly the
classical manner: the head reads the contents of the cells on which it rests, reflects on its state and
follows the rigid instructions of the finite program it is running: accordingly, it writes on the tape,
moves the head one cell to the left or the right or not at all and switches to a new state. Thus, the
classical procedure determines the configuration of the machine at stage $\alpha+1$, given the
configuration at any stage $\alpha$.

We extend the computation into transfinite ordinal time by simply
specifying the behavior of the machine at limit ordinals. When a
classical Turing machine fails to halt, it is usually thought of as
some sort of failure; the result is discarded even though the machine
might have been writing some very interesting information on the tape
(such as all the theorems of mathematics, for example, or the members
of some other computably enumerable set). With infinite time Turing
machines, however, we hope to preserve this information by taking some
kind of limit of the earlier configurations and continuing the
computation transfinitely. Specifically, at any limit ordinal stage
$\lambda$, the head resets to the left-most cell; the machine is placed
in the special {\it limit} state, just another of the finitely many
states; and the values in the cells of the tape are updated by
computing the $\limsup$ of the previous cell values.
\begin{figure}[h]
$$\tapenames\Head101{limit}\Col111\Col010\Col101\Col001\Col011\Dots$$
\caption{The limit configuration}
\end{figure}
With the limit stage configuration thus completely specified, the
machine simply continues computing. If after some amount of time the
{\it halt} state is reached, the machine gives as output whatever is
written on the output tape.

Because there seems to be no need to limit ourselves to finite input and output---the machines have plenty
of time to consult the entire input tape and to write on the entire output tape before halting---the
natural context for these machines is Cantor Space $2^\omega$, the space of infinite binary sequences. For
our purposes here, let's denote this space by $\R$ and refer to its members as real numbers, intending by
this terminology to mean infinite binary sequences. We regard the set of natural numbers $\N$ as a subset
of $\R$ by identifying the number $0$ with the sequence $\<000\cdots>$, the number $1$ with $\<100\cdots>$,
the number $2$ with $\<110\cdots>$, and so on.

Because every program $p$ determines a function---the function sending input $x$ to the output of the
computation of program $p$ on input $x$---the machines provide a model of computation on the reals. We
define that a partial function $f\from\R\to\R$ is {\df infinite time computable} (or {\df supertask
computable}, or for brevity, just {\df computable}, when the infinite time context is understood) when
there is a program $p$ such that $f(x)=y$ if and only if the computation of program $p$ on input $x$ yields
output $y$. A set of reals $A\of\R$ is {\df infinite time decidable} (or {\df supertask decidable} or
again, just {\df decidable}) when its characteristic function, the function with value $1$ for inputs in
$A$ and $0$ for inputs not in $A$, is computable. The set $A$ is infinite time {\df semi-decidable} when
the function of affirmative values $1\restrict A$, that is, the function with domain $A$ and constant value
$1$, is computable. (Thus, the semi-decidable sets correspond in the classical theory to the recursively
enumerable sets, though since here we have sets of reals, we hesitate to describe them as enumerable.)
Since it is an easy matter to change any output value to $1$, the semi-decidable sets are exactly the
domains of the computable functions, just as in the classical theory.

\Theorem. Every supertask computation halts or repeats in countably many steps.

\Proof: Suppose that a supertask computation does not halt by any
countable stage of computation. The point is now that a simple
cofinality argument shows that the complete configuration of the
machine at stage $\omega_1$---the position of the head, the state and
the contents of the cells---must have occurred earlier. For example,
one can find a countable ordinal $\alpha_0$ by which time all of the
cells that have stabilized by $\omega_1$ have stabilized. And then one
can construct a countable increasing sequence of countable ordinals
$\alpha_0<\alpha_1<\cdots$ such that all the cells that change their
value after $\alpha_n$ do so at least once between $\alpha_n$ and
$\alpha_{n+1}$. These ordinals exist because $\omega_1$ is regular and
there are only countably many cells. At the limit stage
$\alpha_\omega=\sup\alpha_n$, which is still a countable ordinal, I
claim that the configuration is the same as at $\omega_1$: since it is
a limit ordinal, the head is on the first cell and in the limit state;
and by construction the contents of each cell are computing the same
$\limsup$ that they compute at $\omega_1$. Since beyond $\alpha_0$ the
only cells that change are the ones that will change unboundedly often,
it follows that limits of this configuration are the very same
configuration again, and the machine is caught in an endlessly
repeating loop. So the proof is complete.\QED

Please observe in this argument that, contrary to the classical
situation, a computation that merely repeats a complete machine
configuration need not be caught in an endlessly repeating loop. After
$\omega$ many repetitions, the limit configuration may allow it to
escape.  One example of this phenomenon would be the machine which does
nothing at all except halt when it is in the {\it limit} state; this
machine repeats its initial configuration many times, yet still halts
at $\omega$.

\section{How powerful are the machines?}

One naturally wants to understand the power of the new machines. The first observation, of course, is that
the classical halting problem for ordinary Turing machines---the question of whether a given program $p$
halts on given input $n$ in finitely many steps---is decidable in $\omega$ many steps by an infinite time
Turing machine. To see this, one programs an infinite time Turing machine to simply simulate the operation
of $p$ on $n$, and if the simulated computation ever halts our algorithm gives the output that yes, indeed,
the computation did halt. Otherwise, the {\it limit} state will be attained, and when this occurs the
machine will that know the simulated computation failed to halt; so it outputs the answer that no, the
computation did not halt.

The power of infinite time Turing machines, though, far transcends the
classical halting problem. The truth is that any question of first
order number theory is supertask decidable. With an infinite time
Turing machine, one could solve the prime pairs conjecture (which
asserts that there are infinitely many primes pairs, pairs of primes
differing by two), for example, and the question of whether there are
infinitely many Fermat primes (primes of the form $2^{2^n}+1$) and so
on: there is a general decision algorithm for any such conjecture.  The
point is that to decide a question of the form $\exists
n\,\varphi(n,x)$, where $n$ ranges over the natural numbers, one can
simply try out all the possible values of $n$ in turn. One either finds
a witness $n$ or else knows at the limit that there is no such witness,
and in this way decides whether $\exists n\,\varphi(n,x)$. Iterating
this idea, one concludes by induction on the complexity of the
statement that any first order number theoretic question is decidable
with only a finite number of limits, that is, before stage $\omega^2$.
In fact, the class of sets that are decidable in time uniformly before
$\omega^2$ is exactly the class of arithmetic sets, the sets of reals
that are definable by a statement using quantifiers over the natural
numbers (see \cite[Theorem 2.6 ]{Hamkins_InfiniteTimeTM}).

\Theorem. Arithmetic truth is infinite time decidable.

One can push this much harder to see that even more complex questions, questions from the lower part of the
projective hierarchy in second order number theory, are supertask decidable. The fact is that any $\Pi^1_1$
set is decidable and more. To prove this, it suffices to consider the most complex $\Pi^1_1$ set, the
well-known set $\WO$, consisting of the reals coding a well-orders of a subset of $\N$. An infinite binary
sequence $x$ codes a relation $\vartriangleleft$ on $\N$ when $i\vartriangleleft j$ if and only if
$x(\<i,j>)=1$, where $\<{\cdot},{\cdot}>$ is the G\"odel pairing function coding pairs of natural numbers
with natural numbers.

\Theorem. The set $\WO$ is infinite time decidable.\label{Count-through}

\Proof: This argument is known as the ``count-through'' argument. We would like to describe a supertask
algorithm which on input $x$ decides whether $x$ codes a well order $\vartriangleleft$ on a subset of $\N$
or not. In $\omega$ many steps, it is easy to check whether $x$ codes a linear order: this amounts merely
to checking that the relation $\vartriangleleft$ coded by $x$ is transitive, irreflexive and connected. For
example, the machine must check that whenever $i\vartriangleleft j$ and $j\vartriangleleft k$ then also
$i\vartriangleleft k$, and all these requirements can be enumerated and checked in $\omega$ many steps.

Next, the algorithm will attempt to find the least element in the field
of the relation $\vartriangleleft$. This can be done by keeping a
current-best-guess on the scratch tape and systematically looking for
better guesses, whenever a new smaller element is found. When such a
better guess is found, it replaces the current guess on the scratch
tape, and a special flag cell is flashed on and then off again. At the
limit, if the flag is on, it means that infinitely often the guess was
changed, and so the relation has an infinite descending sequence. Thus,
in this case the input is definitely not a well order and the
computation can halt with a negative output. Conversely, if the flag is
off, it means that the guess was only changed finitely often, and the
machine has successfully found the $\vartriangleleft$ least element.
The algorithm now proceeds to erase all mention of this element from
the field of the relation $\vartriangleleft$. This produces a new
smaller relation, and the algorithm proceeds to find the least element
of it. In this way, the relation $\vartriangleleft$ is eventually
erased from the bottom as the computation proceeds. If the relation is
not a well order, eventually the algorithm will erase the well founded
initial segment of it, and then discover that there is no least element
remaining, and reject the input. If the relation is a well order, then
the algorithm will eventually erase the entire field, and recognize
that it has done so, and accept the input as a well order. This
completes the proof.\QED

Since $\WO$ is well-known as a complete $\Pi^1_1$ set, we conclude as a corollary that every $\Pi^1_1$ set
is infinite time decidable and hence also, every $\Sigma^1_1$ set is infinite time decidable. But one can't
go much further in the projective hierarchy, because every semi-decidable set has complexity $\Delta^1_2$.
For a finer stratification, let me mention that the arithmetic sets are exactly the sets which can be
decided by an algorithm using a bounded finite number of limits, and the hyperarithmetic sets, the
$\Delta^1_1$ sets, are exactly the sets which can be decided in some bounded recursive ordinal length of
time. Thus, the arithmetic sets are those that can be decided uniformly in time before $\omega^2$, and the
hyperarithmetic sets are exactly those which can be decided uniformly in time before $\omega_1^{ck}$.

Much of the classical computability theory generalizes to the supertask context of infinite time Turing
machines. For example, the s-m-n theorem and the Recursion Theorem go through with virtually identical
proofs. But some other classical results, even very elementary ones, do not generalize. One surprising
result, for example, is the following.

\Theorem. There is a non-computable function whose graph is semi-decidable.\label{Surprise}

This follows from what I have called the Lost Melody Theorem
\cite[Theorem 4.9]{Hamkins_InfiniteTimeTM}, which asserts the existence
of a real $c$ such that $\set{c}$ is decidable, but $c$ is not
writable. Imagine the real $c$ as the melody that you can recognize
when someone sings it, but you cannot sing it on your own. Using such a
lost melody real $c$, one can prove Theorem \ref{Surprise} with the
function $f(x)=c$. Indeed, since this function is constant and the
graph is decidable, the theorem can be strengthened to the assertion
that there is a non-computable constant function whose graph is
decidable. To give some idea of how one proves the Lost Melody Theorem,
let me mention that the real $c$ will be the least real in the G\"odel
constructible universe $L$ hierarchy that codes the ordinal supremum of
the places where all computations on input $0$ have either halted or
repeated. Since this ordinal is above every writable ordinal, the real
$c$ cannot be writable. But the real $c$ codes enough information about
itself so that an infinite time Turing machine can verify that a given
real is $c$ or not.

\section{How long do the computations take?}

One naturally wants to understand how long a supertask computation can take. Therefore, I define an ordinal
$\alpha$ to be {\df clockable} if there is a computation on input $0$ that takes exactly $\alpha$ many
steps to complete (so that the $\alpha^{th}$ step of computation is the act of moving to the halt state).
Such a computation is a clock of sorts, a way to count exactly up to $\alpha$.

It is very easy to see that any finite $n$ is clockable; one can simply have a machine cycle through $n$
states and then halt. The ordinal $\omega$ is clockable, by the machine that halts whenever it sees the
{\it limit} state. And these same ideas show that if $\alpha$ is clockable, then so is $\alpha+n$ and
$\alpha+\omega$. Thus, every ordinal up to $\omega^2$ is clockable. The ordinal $\omega^2$ itself is
clockable: one can recognize it as the first limit of limit ordinals, by flashing a flag on and then off
again every time the {\it limit} state is encountered. The ordinal $\omega^2$ will be first time this flag
is on at a limit stage. Going beyond this, it is easy to see that if $\alpha$ and $\beta$ are clockable, so
are $\alpha+\beta$ and $\alpha\beta$. Undergraduate students might enjoy finding algorithms to clock
specific ordinals, such as $\omega^{\omega^2}$, and I can recommend this as a way to help them understand
the ordinals more deeply.

Most readers will have guessed that the analysis extends much further. In fact, any recursive ordinal is
clockable. This can be seen by optimizing the count-through argument in Theorem \ref{Count-through}.
Specifically, after writing a real coding a recursive ordinal on the tape in $\omega$ many steps, one
proceeds to count through it in an optimized fashion. Rather than merely guessing the least element of the
relation, one guesses the $\omega$ many least elements of the relation (while simultaneously erasing the
previous guesses). In this way, each block of $\omega$ many steps of the algorithm will erase $\omega$ many
elements from the field of the relation.

Some have been surprised that the clockable ordinals extend beyond the recursive ordinals, but in fact they
extend well beyond the recursive ordinals. To see at least the beginnings of this, let me show that the
ordinal $\omega_1^{ck}+\omega$ is clockable, where $\omega_1^{ck}$ is the supremum of the recursive
ordinals. Kleene has proved that there is a recursive relation whose well-founded part has order type
$\omega_1^{ck}$. Consider the supertask algorithm that writes this relation on the tape and then attempts
to count through it. By stage $\omega_1^{ck}$ the ill-founded part will have been reached, but it takes the
algorithm an additional $\omega$ many steps to realize this. So it can halt at stage
$\omega_1^{ck}+\omega$.

One is left to wonder, is $\omega_1^{ck}$ itself clockable? More
generally, {\it Are there gaps in the clockable ordinals?} After all,
if a child can count to twenty-seven, then one might expect the child
also to be able to count to any smaller number, such as
nineteen.\footnote{Friends with children have informed me that such an
expectation is unwarranted; one sometimes can't get the child to stop
at the right time. This reminds me of a time when my younger brother
was in kindergarten, the children all sat in a big circle taking turns
saying the next letter of the alphabet: A, B, C, and so on, around the
circle in the manner of the usual song. After the letter K, the next
child contributed LMNOP, thinking that this was only one letter.} The
question is whether we expect the same to be true for infinite time
Turing machines.

\Theorem. Gaps exist in the clockable ordinals.

\Proof: Consider the algorithm which simulates all programs on input $0$, recording which have halted. When
a stage is found at which no programs halt, then halt. This produces a clockable ordinal above a
non-clockable ordinal, so gaps exist.\QED

The argument can be modified to show that the next gap above any clockable ordinal has size $\omega$. Other
arguments establish that complicated behavior can occur at limits of gaps, because the lengths of the gaps
are unbounded in the clockable ordinals.

\Question. What is the structure of the clockable ordinals?

For example, one might wonder whether the first gap begins at $\omega_1^{ck}$, the supremum of the
recursive ordinals? (It does, since no admissible ordinal is clockable \cite{Hamkins_InfiniteTimeTM}.)

There is another way for infinite time Turing machines to operate as clocks, and this is by counting
through a real coding a well order in the manner of Theorem \ref{Count-through}. To assist with this
analysis, we define that a real is {\df writable} if it is the output of a supertask computation on input
$0$. An ordinal is writable if it is coded by a writable real. It is easy to see that there are no gaps in
the writable ordinals, because if one can write down real coding $\alpha$, it is an easy matter to write
down from this a real coding any particular $\beta<\alpha$. In \cite{Hamkins_InfiniteTimeTM}, Andy Lewis
and I proved that the order types of the clockable and writable ordinals are the same, but the question was
left open as to whether these two classes of ordinals had the same supremum. This was solved by Philip
Welch in \cite{Hamkins_WelchEventual}, allowing Andy Lewis and I to greatly simplify arguments in
\cite{Hamkins_PostProblem}.

\Theorem.(Welch) Every clockable ordinal is writable. The supremum of the writable and clockable ordinals
is the same.

\section{The supertask halting problems}

Any notion of computation naturally provides a corresponding halting problem, the question of whether a
given computation will halt. In the supertask context, we divide the halting problem into two parts, a
boldface and a lightface problem:

$$H=\set{\<p,x> \st \hbox{program $p$ halts on input $x$}}$$
$$h=\set{p \st \hbox{program $p$ halts on input $0$}}$$

In the classical theory, of course, these two sets are Turing equivalent, but here the situation is
different. Nevertheless, for undecidability the classical arguments do directly generalize.

\Theorem. The halting problems $h$ and $H$ are semi-decidable but not decidable.

For semi-decidability, the point is that given a program $p$ and input $x$ (or input $0$), one can simply
simulate $p$ on $x$ to see if it halts. If it does, output the answer that yes, it halted; otherwise, keep
simulating. For undecidability, in the case of $H$ one can use the classical diagonalization argument; for
the lightface halting problem $h$, one appeals to the Recursion Theorem, just as in the classical theory.

\section{Oracles}

There are two natural types of oracles to use in the infinite time Turing machine context. On the one hand,
one can use an individual real as an oracle just as one does in the classical context, by simply adding an
oracle tape containing this real, and allowing the machine to access this tape during the computation. This
corresponds exactly to adding an extra input tape and thinking of the oracle real as a fixed additional
input.

But this is ultimately not the right type of oracle to consider. Rather, an oracle is more properly the
same type of object as one that might be decidable or semi-decidable, namely, a {\it set} of reals, not an
individual real. Since such a set could be uncountable, we can't expect to be able to write out the entire
contents of the oracle on an extra tape. Rather, we provide an oracle model of relative computability by
which the machine can make arbitrary membership queries of the oracle. Specifically, for a fixed oracle set
of reals $A$, we equip an infinite time Turing machine with an initially blank oracle tape on which the
machine can read or write. By attempting to switch to a special {\it query} state, the machine receives the
answer (by moving actually to the {\it yes} or {\it no} state) as to whether the real currently written on
the oracle tape is in $A$ or not. In this way, the machine is able to ask, of any real $x$ that it is
capable of producing, whether $x\in A$ or not. This model of oracle computation has proven robust, and it
closely follows the well-known definition of $L[A]$ in set theory, the constructible universe relative to
the predicate $A$, in which at any given stage in the construction one is allowed to apply the predicate
only to previously constructed objects.

From the notion of oracle computation, one can of course define a notion of relative computability.
Specifically, the set $A$ is computable from $B$, written {\df $A\ileq B$}, if and only if $A$ is supertask
decidable using oracle $B$. One then also defines $A\iequiv B$ if and only if $A\ileq B$ and $B\ileq A$,
and this is the equivalence relation of the infinite time Turing degrees. The strict version $A\ilt B$
holds if and only if $A\ileq B$ and $A\not\iequiv B$.

\section{Supertask Jump Operators}

The two halting problems give rise of course to two jump operators. Specifically, for any set $A$ we have
the boldface and lightface jumps:
$$A^\Jump=H^A=\set{\<p,x>\st\hbox{program $p$ halts on input $x$ with oracle $A$}}$$
$$A^\jump=A\oplus h^A=A\oplus\set{p\st\hbox{program $p$ halts on input $0$ with oracle $A$}}$$
We include the factor $A$ explicitly in $A^\jump$, because in general $A$ may not be computable from $h^A$.
Indeed, there are some sets $A$ that are not computable from any real at all.

\Theorem Jump Theorem. For any set, $A\ilt A^\jump\ilt A^\Jump$.

To prove this theorem, one first observes that $A\ileq A^\jump\ileq A^\Jump$, since $A$ is explicitly
computable from $A^\jump$ and $A^\jump$ is merely the $0^{th}$ slice of $A^\Jump$. Secondly, one knows that
$A\ilt A^\jump$ because the undecidability of the relativized halting problem means that $h^A$ is not
computable from $A$. The nontrivial aspect of this theorem is the assertion that $A^\jump\ilt A^\Jump$.
This assertion is what separates the two jump operators, and is the reason that we know the two halting
problems $h\iequiv 0^\jump$ and $H\iequiv 0^\Jump$ are not equivalent. This follows from the more specific
result that the set $A^\Jump$ is not computable from $A\oplus z$ for any real $z$. In particular, $0^\Jump$
is not computable from any real. In fact the boldface jump $\Jump$ jumps much higher than the lightface
jump $\jump$, and absorbs many iterates of the weaker jump, since $A^{\jump\Jump}\iequiv A^\Jump$; indeed,
for any ordinal $\alpha$ which is $A^\Jump$-writable, $A^{\jump^{(\alpha)}\Jump}\iequiv A^\Jump$ (see
\cite{Hamkins_InfiniteTimeTM}).

\section{Post's Problem for Supertasks}

Post's problem is the question in classical computability theory of
whether there are any non-decidable semi-decidable degrees strictly
below the halting problem, or equivalently, whether there are any
intermediate semi-decidable degrees between $0$ and the Turing jump
$0'$. This question has a natural supertask analogue:

\Question Supertask Post's Problem. Are there any intermediate
semi-decidable supertask degrees between $0$ and the supertask jump
$0^\jump$?

The answer is delicately mixed. On the one hand, in the context of degrees in the real numbers, we have a
negative answer. This contrasts sharply with the classical theory.

\Theorem. There are no reals $z$ such that $0\ilt z\ilt 0^\jump$.

\Proof: Suppose that $0\ileq z\ileq 0^\jump$. So $z$ is the output of program $p$ using $0^\jump$ as an
oracle. Consider the algorithm which computes approximations to $0^\jump$, and uses program $p$ with these
approximations in an attempt to produce $z$. If one of the proper approximations to $0^\jump$ can
successfully produce $z$, then $z$ is writable and $0\iequiv z$. Conversely, if none of the proper
approximations can produce $z$, then on input $z$ we can recognize $0^\jump$ as the {\it true
approximation}, the first approximation able to produce $z$. So $z\iequiv 0^\jump$.\QED

On the other hand, when it comes to {\it sets} of reals, we have an affirmative answer.

\Theorem. There are semi-decidable sets of reals $A$ with $0\ilt A\ilt
0^\jump$. Indeed, there are incomparable semi-decidable sets $A\perp
B$.\label{PostProblem}

Please consult \cite{Hamkins_InfiniteTimeTM} for the proof. Let me
mention here, though, that the basic idea of the argument is to
generalize the Friedburg-Munchnik priority argument to the supertask
context, much as Sacks' did for $\alpha$-recursion theory. Building $A$
and $B$ in stages, we attempt to meet the requirements
$$\varphi_p^B\not=A\qquad{\ \rm and\ }\qquad\varphi_p^A\not=B$$
by adding writable reals to $A$ and $B$ that have not yet appeared on
the higher priority computations. One technical fact to make this idea
work is that for any clockable ordinal $\alpha$, there are many
writable reals not appearing during the course of any supertask
computation of length $\alpha$. Thus, we can find a supply of new
writable reals to add to $A$ and $B$ in order to satisfy the later
requirements, without injuring the witnessing computations of earlier
higher-priority requirements.

\section{Other Models of Infinitary Computation}

Let me briefly compare the infinite time Turing machine model of
supertask computation with some other well-known models.

The Blum-Shub-Smale machines (see \cite{Hamkins_BlumShubSmale}) were
the original inspiration for infinite time Turing
machines.\footnote{Jeff Kidder and I heard Lenore Blum's lectures for
the Berkeley Logic Colloquium in 1989, and had the idea to generalize
the Turing machine concept in a different direction: to infinite time
rather than infinite precision.} Programs and computations for BSS
machines are finite, but the basic units of computation are full
precision real numbers. They are in essence finite state register
machines, where the registers each hold a real number. The primary
purpose of introducing the BSS machines was to provide a theoretical
foundation for analyzing computational algorithms using the concepts of
real analysis rather than arithmetic. The machines allow one to analyze
the dynamical features, for example, of actual algorithms in numerical
analysis, such as Newton's method, and illuminate questions of
stability and convergence for such algorithms. The classical approach
to these problems, using the Turing machine model with ever greater
decimal approximations, forces one into the realm of finite
combinatorics, where one becomes lost in a jumble of discrete
approximation error analysis, when one would rather fly smoothly above
it in the heaven of differential equations.

In another direction, the theory of higher recursion provides a model
of infinitary computability by setting a very general theoretical
context for recursion on infinite objects, and one should expect many
parallels between it and the theory of infinite time Turing machines.
The anonymous referee of \cite{Hamkins_PostProblem} and Philip Welch
have pointed out, for example, that the infinitary priority argument
\cite[Theorem 4.1]{Hamkins_PostProblem}, stated as Theorem
\ref{PostProblem} above, parallels Sacks' version of the
Friedburg-Munchnik proof for $\alpha$-recursion \cite{Hamkins_Sacks},
specifically when $\alpha$ is $\lambda$, the supremum of the clockable
ordinals. One can identify the writable reals in our argument with the
ordinal stages at which they appear and get Sacks' sets, and
conversely, Sacks' could have written out codes for those stages and
gotten our sets.  This identification reveals that the $\ileq$-degree
structure of sets of writable reals below $0^\jump$ is exactly that of
the $\lambda$-degrees. Accordingly, one can obtain not only the answer
to Post's problem, but all the theorems from $\lambda$-recursion theory
for this class of degrees, such as the Shore Density Theorem, etc., for
free. It will be very interesting to see if these ideas will allow one
to prove the theorems in the general case of all degrees.

Lastly, let me mention quantum Turing machines, if only because I am
often asked about them in connection with infinite time Turing
machines. Quantum Turing machines are like classical Turing machines,
except that the configuration of the machine at any given stage is a
superposition of classical configurations; the different components of
these superpositions, like the wave functions of quantum mechanics, may
constructively or destructively interfere with one another as the
computation proceeds. By means of clever quantum algorithms, one can
effectively carry out parallel computation in these different
components, constructively interfering their output to assemble the
information into a final answer. In this way, quantum Turing machines
allow for an exponential increase in the speed of computation of many
important functions. But because quantum Turing machines, at the end of
the day, are simulable by classical Turing machines, they do not
introduce new decidable sets or new computable functions. And so while
quantum Turing machines are without a doubt extremely important in
matters of computational feasibility, they do not really provide a
model of infinitary computability. Infinite time Turing machines are
simply much more powerful than quantum Turing machines.

\section{Questions for the Future}

I close this article by asking the open-ended question:

\Question. What is the structure of infinite time Turing degrees? To what extent do its properties mirror
or differ from the classical structure?

This question really stands for the dozens of specific open questions that one might ask: does the Sacks
Density Theorem, for example, hold in the supertask context for arbitrary sets of reals? The field is wide
open.

\articleend